\documentclass{article}

\usepackage{a4wide}

\usepackage{amsfonts}
\usepackage{amssymb}
\usepackage{graphicx}
\usepackage[centertags]{amsmath}
\usepackage{amsthm}
\usepackage{newlfont}
\usepackage{lscape}
\usepackage{setspace}
\usepackage[hmargin={3.2cm}]{geometry}
\usepackage{subfig}

\begin{document}

\title{Optimal control of a dengue epidemic model
with vaccination\thanks{This is a preprint
of a paper accepted for presentation at ICNAAM 2011,
Halkidiki, Greece, 19-25 September 2011,
and to appear in \emph{AIP Conference Proceedings}, volume 1389.}}

\author{Helena Sofia Rodrigues\\
School of Business Studies, Viana do Castelo Polytechnic Institute, Portugal\\
\texttt{sofiarodrigues@esce.ipvc.pt}
\and
M. Teresa T. Monteiro\\
Department of Production and Systems, University of Minho, Portugal\\
\texttt{tm@dps.uminho.pt}
\and
Delfim F. M. Torres\\
Center for Research and Development in Mathematics and Applications\\
Department of Mathematics, University of Aveiro, Portugal\\
\texttt{delfim@ua.pt}}

\date{}

\maketitle

% -------------------------------------------

\begin{abstract}
We present a SIR+ASI epidemic model to describe
the interaction between human and dengue fever mosquito populations.
A control strategy in the form of vaccination,
to decrease the number of infected individuals, is used.
An optimal control approach is applied in order
to find the best way to fight the disease.

\bigskip

\noindent \textbf{Keywords:} optimal control,
Pontryagin maximum principle, dengue, vaccination.

\smallskip

\noindent \textbf{MSC 2010:} 49K15, 92D30.

\smallskip

\noindent \textbf{PACS:} 87.23.Cc, 87.55.de
\end{abstract}

% -------------------------------------------

\section{Introduction}
\label{introduction}

Dengue fever is a vector borne disease, which has become an
increasingly public health problem that carries a huge financial
burden to the governments. Currently, the only way of
controlling the disease is to minimize the vector population.
Dengue vaccine for effective prevention and long
term control under development, is expected to be the solution.
Dengue transcends international borders
and is emerging rapidly as a consequence of globalization
and climate changes. It is a disease of great complexity,
due to the interactions between humans, mosquitoes,
and various virus serotypes as well as
efficient vector survival strategies. The four serotypes, known as
DEN1 to DEN4, constitute a complex of flaviviridae transmitted by
\emph{Aedes} mosquitos, specially \emph{Aedes Aegypti}. Infection by
any of the four serotypes induces lifelong immunity against
reinfection by the same type, but only partial and temporary
protection against the others. Sequential infection by different
serotypes could lead to a more severe dengue episode: dengue
hemorrhagic fever (DHF).

Vector control remains the only available strategy against dengue.
Despite integrated vector control with community participation,
along with active disease surveillance and insecticides,
there are only a few examples of successful dengue prevention
and control on a national scale \cite{Cattand2006}.
To make matters worse, the levels of resistance of \emph{Aedes Aegypti} to
insecticides has increased, which imply shorter intervals between
treatments, and only few insecticide products are available in the
market because of high costs for development and registration and low returns.

For long time, the evaluation of global dengue disease burden
was limited and the stakeholders considered the
potential market for the dengue vaccine to be small. By the end of
20th century, with the increase in dengue infections as well as
the prevalence of all four circulating serotypes,
faster development of a vaccine became a serious concern \cite{Murrell2011}.
It is agreed that a vaccination program not only protects directly the individual,
but also indirectly the population, which is called herd immunity.
As a consequence of vaccination, the occurrence of epidemics would decrease
relieving health facilities. However, constructing a successful vaccine for dengue
has been challenging. Not only is the knowledge of disease pathogenesis
insufficient, but also the vaccine must protect against
all serotypes so that the level of DHF doesn't increase.

% ------------------------------------------

\section{Optimal control of the epidemiological model}
\label{epidemiological model}

Two types of population were considered: hosts and vectors.
The hosts (humans) are divided into three complementary classes:
susceptible, $S_h (t)$, individuals who can contract the disease;
infected, $I_h(t)$, individuals capable of transmitting the
disease to others; and resistant, $R_h (t)$, individuals who have
acquired immunity at time $t$. The total number of hosts is
constant, which means that $N_h=S_h(t)+I_h(t)+R_h(t)$. Similarly, the
model has also three compartments for the vectors (mosquitoes):
$A_m(t)$, which represents the aquatic phase of the mosquito
(including egg, pupae and larvae) and the adult phase of the
mosquito, with $S_m(t)$ and $I_m(t)$, susceptible and infected,
respectively. It is also assumed that $N_m=S_m(t)+I_m(t)$.
The model is described by an initial value problem with a system
of six differential equations:
\begin{equation}
\label{ode}
\begin{cases}
\frac{dS_h}{dt} = \mu_h N_h- \left(B\beta_{mh}\frac{I_m}{N_h}+\mu_h+u\right)S_h+\sigma u R_h\\
\frac{dI_h}{dt} = B\beta_{mh} \frac{I_m}{N_h} S_h -(\eta_h+\mu_h) I_h\\
\frac{dR_h}{dt} = \eta_h I_h + u S_h - \left(\sigma u+\mu_h\right) R_h\\
\frac{dA_m}{dt} = \varphi \left(1-\frac{A_m}{k N_h}\right) (S_m+I_m) - \left(\eta_A+\mu_A\right) A_m\\
\frac{dS_m}{dt} = \eta_A A_m - \left(B \beta_{hm}\frac{I_h}{N_h}+\mu_m\right) S_m\\
\frac{dI_m}{dt} = B \beta_{hm}\frac{I_h}{N_h}S_m -\mu_m I_m .
\end{cases}
\end{equation}
The recruitment rate of human population is noted by $\mu_hN_h$.
The natural death rate for humans and mosquitoes, aquatic and
adult phase, is described by the parameters $\mu_h$, $\mu_m$ and $\mu_A$,
respectively. We assume that $B$ is the average daily
biting (per day) of the mosquito whereas $\beta_{mh}$ and
$\beta_{hm}$ are related to the transmission probability (per
bite) from infected mosquitoes to humans and vice versa.
By $\varphi$ we denote the number of eggs at each deposit per capita (per day).
The recovery rate of the human population is denoted by $\eta_h$. The
maturation rate from larvae to adult (per day) is denoted by
$\eta_{A}$. The vaccine coverage of the susceptible is represented
by $u$ (the control variable). The factor $\sigma$ represents the
level of inefficacy of the vaccine: for $\sigma=0$ the vaccine is
perfectly effective, while $\sigma=1$ means that the vaccine has no
effect at all. The main aim of this work is to study
the optimal vaccination strategy considering both
the costs of treatment of infected individuals and the costs
of vaccination. So, the objective functional is
\begin{equation}
\label{functional}
\text{ minimize } J[u]=\int_{0}^{t_f}\left[\gamma_I I_h(t)^2
+\gamma_V u(t)^2\right]dt,
\end{equation}
where $\gamma_I$ and $\gamma_V$ are positive constants
representing the weights of the costs
of treatment of infected and vaccination, respectively.
Let $\lambda_i(t)$, with $i=1,\ldots,6$, be the co-state variables.
The Hamiltonian for the present optimal control problem is given by
\begin{equation}
\label{hamiltonian}
\begin{split}
H &= \lambda_1\left[\mu_h N_h- \left(B\beta_{mh} \frac{I_m}{N_h}
+ \mu_h+u\right)S_h +\sigma u R_h\right]
+\lambda_2\left[B\beta_{mh} \frac{I_m}{N_h} S_h
-\left(\eta_h+\mu_h\right) I_h\right]\\
&+ \lambda_3 \left[\eta_h I_h +uS_h- \left(\sigma u+\mu_h\right) R_h\right]
+\lambda_4\left[\varphi \left(1-\frac{A_m}{k N_h}\right)(S_m+I_m)-\left(\eta_A+\mu_A\right)A_m\right]\\
&+ \lambda_5\left[\eta_A A_m -\left(B \beta_{hm}\frac{I_h}{N_h}+\mu_m\right) S_m\right]
+\lambda_6\left[B \beta_{hm}\frac{I_h}{N_h}S_m -\mu_m I_m\right] + \gamma_I I_h^2+\gamma_V u^2.
\end{split}
\end{equation}
By the Pontryagin maximum principle \cite{book:Pont}, the optimal control $u^{*}$
should be the one that minimizes, at each instant $t$, the Hamiltonian given
by \eqref{hamiltonian}, that is,
$$
H\left(x^*(t),\lambda^*(t),u^*(t)\right)
=\min_{u \in [0,1]} H\left(x^*(t),\lambda^*(t),u\right).
$$
The optimal control, derived from the stationary condition
$\frac{\partial H}{\partial u}= 0$ and considering
$0\leq u\leq 1$, is given by
\begin{equation*}
u^*=\min\left\{1,\max\left\{0,\frac{\left(\lambda_1
-\lambda_3\right)\left(S_h-\sigma R_h\right)}{2\gamma_V}\right\}\right\}.
\end{equation*}
Substituting the optimal control $u^*$ into the state system \eqref{ode}
and the adjoint system $\lambda^{'}_i(t)=-\frac{\partial H}{\partial x_i}$,
\textrm{i.e.},
\begin{equation*}
\begin{cases}
\frac{d\lambda_1}{dt}=(\lambda_1-\lambda_2) \left(B \beta_{mh} \frac{I_m}{N_h}\right)
+\lambda_1 \mu_h+(\lambda_1-\lambda_3)u\\
\frac{d\lambda_2}{dt}=-2\gamma_I I_h+\lambda_2(\eta_h+\mu_h)
-\lambda_3\eta_h+(\lambda_5-\lambda_6)\left(B\beta_{hm}\frac{S_m}{N_h}\right)\\
\frac{d\lambda_3}{dt}=-\lambda_1\sigma u+\lambda_3(\mu_h+\sigma u)\\
\frac{d\lambda_4}{dt}=\lambda_4 \varphi \frac{S_m+I_m}{k N_h}
+\lambda_4(\eta_A+\mu_A)-\lambda_5\eta_A\\
\frac{d\lambda_5}{dt}=-\lambda_4 \varphi \left(1-\frac{A_m}{k N_h}\right)
+(\lambda_5-\lambda_6) B\beta_{hm} \frac{I_h}{N_h}+\lambda_5\mu_m\\
\frac{d\lambda_6}{dt}=(\lambda_1-\lambda_2) \left(B\beta_{mh}\frac{S_h}{N_h}\right)
-\lambda_4 \varphi \left(1-\frac{A_m}{k N_h}\right) + \lambda_6 \mu_m,
\end{cases}
\end{equation*}
we obtain the corresponding $x^*$ and $\lambda_i^*$, $i=1,\ldots 6$,
with the help of the transversality conditions
$\lambda_i^*(t_f)=0$, $i=1,\ldots,6$
(see \cite{book:Pont} for details).

% ------------------------------------------

\section{Numerical simulation and discussion}
\label{numerical}

\begin{figure}[!ht]
\begin{minipage}[b]{.5\linewidth}
\centering
\includegraphics[scale=.25]{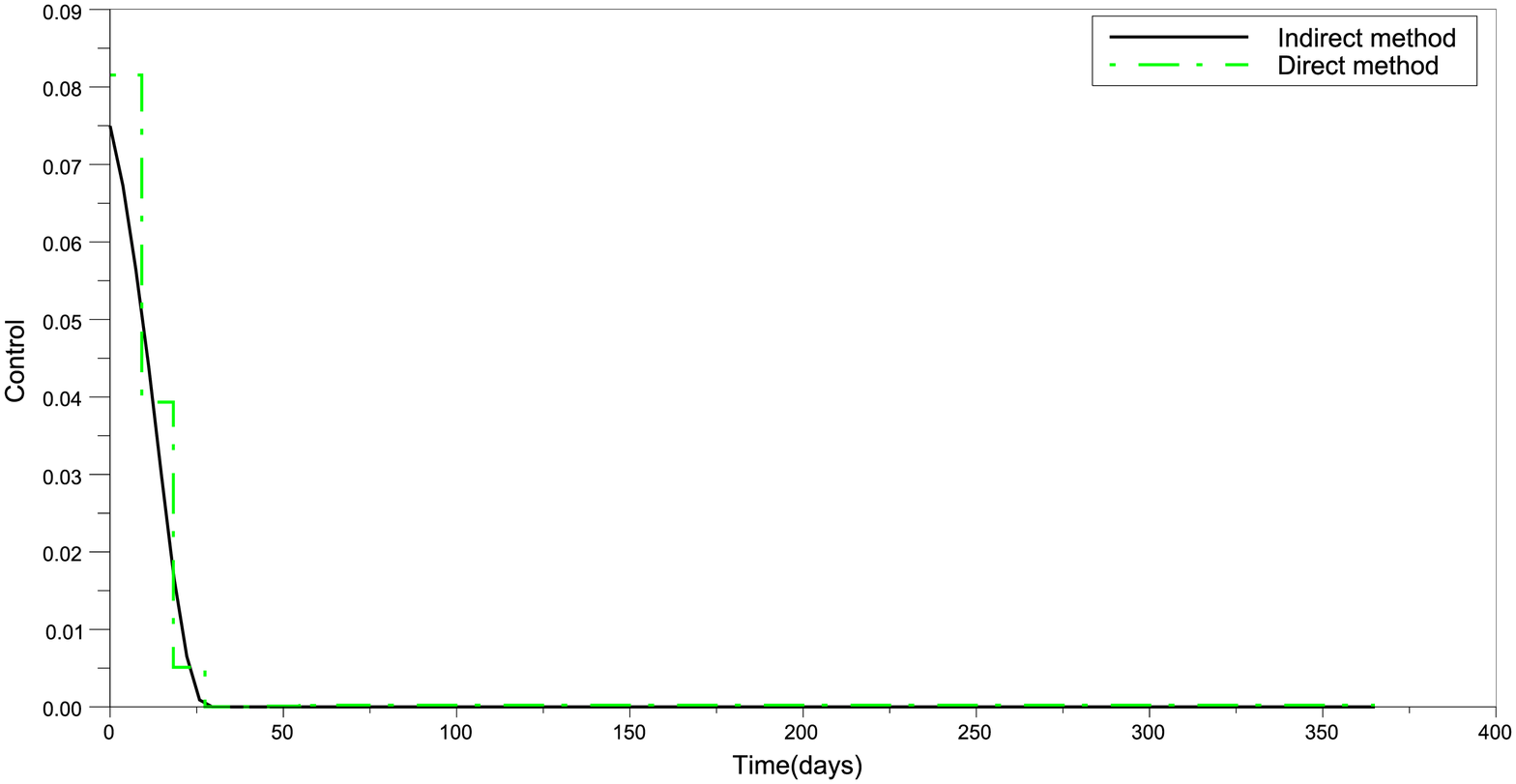}
\caption{Optimal control with direct and indirect approaches.} \label{Fig1}
\end{minipage}
\hspace{.3cm}
\begin{minipage}[b]{.5\linewidth}
\centering
\includegraphics[scale=.25]{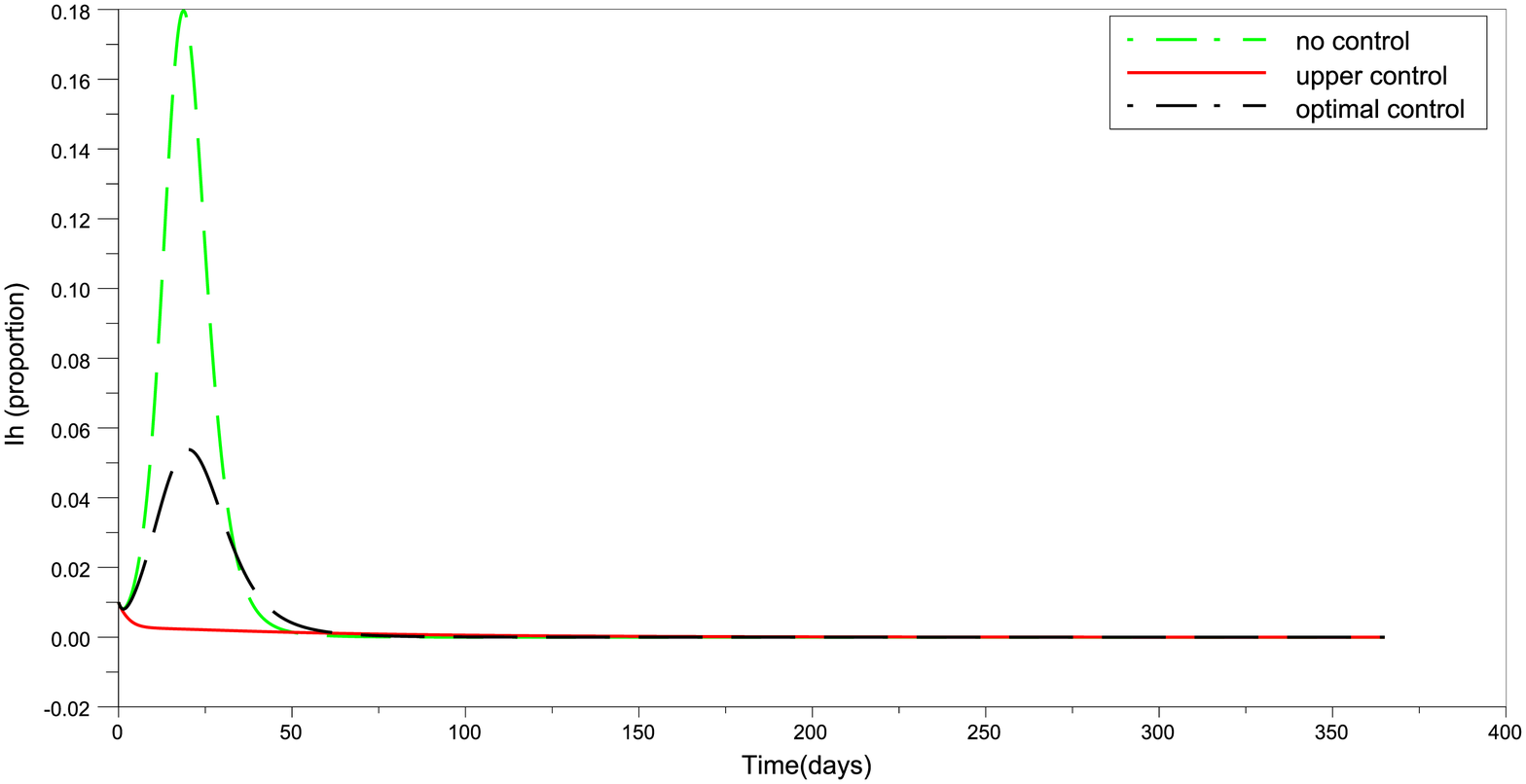}
\caption{Infected humans using different levels of control.} \label{fig1}
\end{minipage}
\end{figure}

The simulations were carried out using the following values:
$N_h = 480000$, $B = 0.5$, $\beta_{mh} = 0.3$, $\beta_{hm} = 0.3$,
$\mu_h = 1/(71 \times365)$, $\eta_h = 1/3$, $\mu_m = 1/10$, $k=3$,
$m=3$, $N_m=m\times N_h$, $\varphi = 6$, and $t_f=365$ days.
It was considered that the vaccine is imperfect with a level of inefficacy of
$\sigma=0.15$. The initial conditions for the ordinary
differential system were: $S_h(0) = N_h-216$, $I_h(0)=216$,
$R_h(0)=0$, $A_m=k*N_h$, $S_m(0)=N_m$ and $I_m(0)=0$.
The optimal control problem was solved using two methods: direct
and indirect. For an introduction to direct and indirect
methods in optimal control we refer the reader to \cite{MR2589623,MR2224013}.
The direct method uses the optimal functional
\eqref{functional} and the state system \eqref{ode} and was solved by
\textsf{DOTcvpSB} \cite{Dotcvp}. It is a toolbox implemented
in \textsf{MatLab}, which uses an ensemble of numerical methods
for solving continuous and mixed-integer dynamic optimization problems.
The indirect method we used is an iterative method with a Runge--Kutta scheme,
solved through \emph{ode45} of \textsf{MatLab}. The state system with an
initial guess is solved forward in time and then the adjoint
system with the transversality conditions is solved backward in
time. The controls are updated at the end of each iteration (see
\cite{Lenhart2007} for more details).
Figure~\ref{Fig1} shows the optimal control obtained by
the two different approaches. They both seem to
have the same behavior.
% ------------------------------------------
% ------------------------------------------
\begin{table}
\begin{tabular}{lccc}
\hline \\Method & optimal control & no control ($u \equiv 0$)& upper control ($u \equiv 1$)\\
\hline Direct (DOTcvpSB) & 0.146675 & 0.674555 & 364.940488\\
Indirect (backward-forward) & 0.113137 & 0.357285 & 365.00046\\
\hline
\end{tabular}
\caption{Values of the cost functional \eqref{functional}} \label{resultados}
\end{table}
% ------------------------------------------
Table~\ref{resultados} shows the costs obtained by the two
methods, in three situations: optimal control, no control
($u(t) \equiv 0$) and upper control ($u(t) \equiv 1$).
Figure~\ref{fig1} shows the number of infected humans when
different controls are considered. It is possible to see that the
upper control, which means that everyone is vaccinated, implies
that just a few individuals were infected, allowing eradication of the
disease. Although the optimal control, in the sense of objective
\eqref{functional}, allows the occurrence of an outbreak,
the number of infected individuals is much lower when
compared with a situation where no one is vaccinated.
Also, the costs are very low when compared with the upper control case.

% ------------------------------------------

\section{Conclusions}
\label{conclusions}

Dengue is an infectious tropical disease difficult to prevent and manage.
Researchers agree that the development of a vaccine
for dengue is a question of high priority.
In the present study we show how a vaccine results
in saving lives and at the same time in a
reduction of the budget related with the disease.
As future work we intend to study the interaction of a dengue
vaccine with other kinds of control already investigated in the literature,
such as insecticide and educational campaigns
\cite{Rodrigues2010b,Rodrigues2011}.

% ------------------------------------------

\section*{Acknowledgments}

Work partially supported by the Portuguese Foundation for Science
and Technology (FCT) through the Ph.D. grant SFRH/BD/33384/2008
(Rodrigues) and the R\&D units Algoritmi (Monteiro) and CIDMA (Torres).

% ------------------------------------------

% ------------------------------------------

\end{document}